\documentclass[10pt]{article}
\usepackage[dvips]{graphicx}
\usepackage{amssymb}
\usepackage{amsmath}
\usepackage{caption}
\usepackage{hyperref}
\usepackage{url} 
\usepackage{arcs}
\usepackage{xcolor}
\usepackage[scr=boondoxo,scrscaled=1.05]{mathalfa}
\usepackage{mathptmx} 
\usepackage{multirow}

\newtheorem{Remark}{Remark}

\newcommand{\hybridhyper}{\mathscr{H}_{n}^{\lambda}}

\newcommand{\discreteinner}[1]{\left\langle#1\right\rangle_M}

\newcommand{\softoperatordef}[1]{\mathscr{S}_{k}\left({#1}\right)}

\newenvironment{AMS}{\small\bf 2010 AMS subject classification: }{}

\begin{document}

\title{\bf
Numerical cubature and hyperinterpolation over Spherical Polygons
\thanks{Work partially
supported by the
DOR funds 
of the University of Padova, and by the INdAM-GNCS grant {\it{Kernel and polynomial methods for approximation and integration: theory and software for applications}}.
This research has been accomplished within the RITA {\it{Research ITalian network on Approximation}}, the SIMAI Activity Group ANA\&A and the UMI Group TAA {\it{Approximation Theory and Applications}}.
}}

\author{A. Sommariva\\
\small{University of Padova, Italy}}


\maketitle

\footnotetext[1]{Corresponding author: alvise@math.unipd.it}

\begin{abstract}
The purpose of this work is to introduce a strategy for determining the  nodes and weights of a low-cardinality positive cubature formula nearly exact for polynomials of a given degree over spherical polygons.
In the numerical section we report the results about numerical cubature over a spherical polygon $\cal P$ approximating Australia and reconstruction of functions over such $\cal P$, also affected by perturbations, via hyperinterpolation and some of its variants. The open-source Matlab software used in the numerical tests is available at the author's homepage.
\end{abstract}

\vskip0.2cm
\noindent
\begin{AMS}
{\rm Primary 41A10,  42A10, 65D05, 65D32.}
\end{AMS}

\section{Introduction}

In \cite{SV21A} the authors introduced a numerical code for the computation of cubature formula on spherical triangles of the unit-sphere ${\mathbb{S}}_2$, nearly exact for polynomials $p_n \in {\mathbb{P}}_n({\mathbb{S}}_2)$, i.e. numerically achieving a given total polynomial degree $n$. 
The algorithm was based on subperiodic trigonometric gaussian quadrature for planar elliptical sectors and on Caratheodory-Tchakaloff quadrature compression via NNLS. The final rule had internal nodes, positive weights and cardinality at
 most equal to $(n+1)^2$, that is the dimension of the vector space  ${\mathbb{P}}_n({\mathbb{S}}_2)$. 

Later in \cite{SV21B}, they move their attention to the computation of an orthogonal polynomial basis on spherical triangles, via the formula described above, and to the construction of the corresponding weighted orthogonal projection 
(hyperinterpolation) of a function sampled at the cubature nodes. 

Having in mind the studies of real world problems, e.g. coming from geo-mathematics, we move here to the more general setting of spherical polygons $\cal P$ of ${\mathbb{S}}_2$. 

We first show how to implement such a cubature rule nearly exact for polynomials in $ {\mathbb{P}}_n({\mathbb{S}}_2)$, with internal nodes, positive weights and cardinality at most equal to $(n+1)^2$. Next we investigate its qualities when applied to certain spherical polygons, reporting cardinalities, cputimes as well as cubature errors relatively to certain integrands with different regularity.

Later we give a brief introduction to the classical hyperinterpolation and to some of its variants useful in case of noisy data, finally testing the numerical reconstruction of functions also in the case of perturbations given by impulsive and gaussian noise.

\section{Numerical cubature on spherical polygons}
We start our study introducing some definitions that are useful in the paper.

A {\sl{great circle}} is the intersection of the unit-sphere ${\mathbb{S}}^2=\{(x,y,z): x^2+y^2+z^2=1\}$ with a plane passing through the origin.

Let $V_1, \ldots,V_L$ be distinct points of ${\mathbb{S}}^2$ and set $V_0=V_L$. A {\sl{spherical polygon}} with vertices $\{V_k\}_{k=1,\ldots,L}$ is the region ${{\cal{P}}} \subset {\mathbb{S}}^2$ whose boundary $\delta {{\cal{P}}}$  is determined by the geodesic arcs $\{\gamma_k\}_{k=0,\ldots,L}$, where each edge $\gamma_k$ is the portion of the great circle joining $V_k$ with $V_{k+1}$, and is oriented counterclockwise.

In this paper we suppose that the spherical polygon ${{\cal{P}}}$ is contained in a cap whose {\sl{polar angle}} is strictly inferior than $\pi$. We notice that, for the purpose of determining a cubature rule over ${{\cal{P}}}$, if the spherical polygon has not this feature then it can always be subdivided in spherical polygonal regions ${{\cal{P}}}_k$, each one contained in a cap whose {\sl{polar angle}} is strictly inferior than $\pi$, whose interior do not overlap and whose union is ${{\cal{P}}}$. Next, one can apply the procedure that we will explain below on each ${{\cal{P}}}_k$, obtaining a {\sl{composite}} rule on ${{\cal{P}}}$.

In this section we show how to obtain a cubature rule of PI-type over such a spherical polygon ${{\cal{P}}}$, i.e. a formula that has {\sl{positive weights}} $\{w_k\}_{k=1,\ldots,M}$ and {\sl{internal nodes}} $\{Q_k\}_{k=1,\ldots,M}$, and that is nearly exact in the space ${\mathbb{P}}_n({\mathbb{S}}_2)$ of polynomials on ${\mathbb{S}}_2$ of total degree at most $n$.

It can be proven that if ${{\cal{P}}} \subset {\mathbb{S}}_2$ is a spherical polygon, each of whose edges have length strictly less than $\pi$, and ${{\cal{P}}}$ does not contain a great circle, then ${{\cal{P}}}$ has a spherical triangulation \cite{OR}. Consequently, since ${{\cal{P}}}$ is contained in a cap whose {\sl{polar angle}} is strictly inferior than $\pi$, we can conclude that there is a spherical triangulation of ${{\cal{P}}}$. 

{\vspace{0.25cm}}

The latter can be obtained as follows:
\begin{enumerate}
\item denoting by $C_{{\cal{P}}}$ the centroid of ${{\cal{P}}}$, we determine the {\sl{gnomonic projection ${\tilde{{\cal{P}}}}$ of ${{\cal{P}}}$ on the plane $\pi_{C_{{\cal{P}}}}$ tangent in $C_{{\cal{P}}}$}} to the unit-sphere;
\item {\sl{compute a triangulation over the planar polygon ${\tilde{{\cal{P}}}}$}} (e.g. by Matlab built-in environment {\tt{polyshape}} that determines a minimal triangulation), i.e. ${\tilde{{\cal{P}}}}=\cup_{i=1}^m{\tilde{{{\cal{T}}}}}_i$ where ${\tilde{{\cal{T}}}}_i \subset \pi_{C_{{\cal{P}}}}$  are planar triangles whose interiors do not overlap, i.e. if $j \neq k$ then ${\mbox{int}}({\tilde{{\cal{T}}}}_j) \cap {\mbox{int}}({\tilde{{\cal{T}}}}_k)=\emptyset$;
\item {\sl{map back to the sphere}}, by means of the inverse of the gnomonic projection, each planar triangle ${\tilde{{\cal{T}}}}_k$ into a spherical triangle ${\cal{T}}_k$, obtaining the required spherical triangulation.
\end{enumerate}

Since $\{{\cal{T}}_k\}_{k=1,\ldots,m}$ is a spherical triangulation of ${{\cal{P}}}$, if we determine a cubature rule of PI-type with ADE equal to $n$ over each spherical triangle ${\cal{T}}_k$ then by the addivity of integration we have immediatly a rule with the same features also on ${{\cal{P}}}$. 

\begin{figure}[!ht]
\centering
{\includegraphics[scale=0.4,clip]{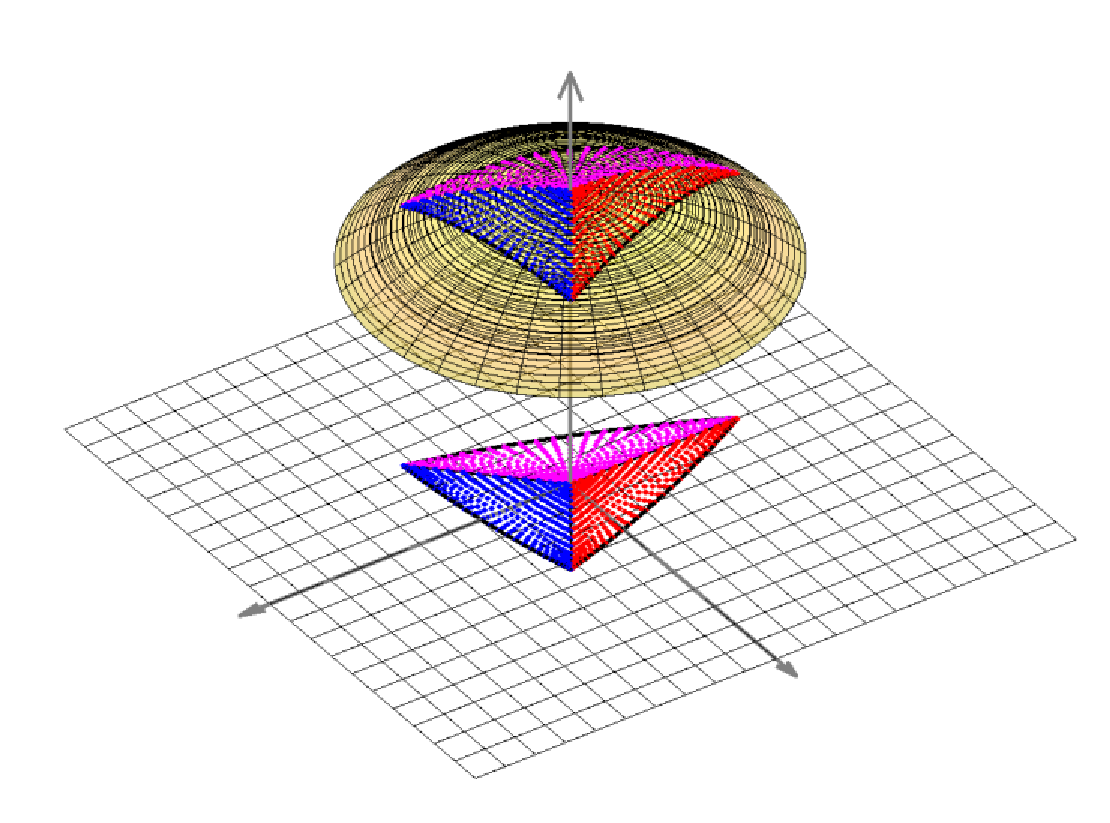}}
 \caption{cubature nodes (and weights) on a spherical triangle lifted from the projected elliptical triangle.}
 \label{fig1}
 \end{figure}

For this reason, we concentrate our efforts in determining a rule over a spherical triangle ${{\cal{T}}}$ contained in a cap whose {\sl{polar angle}} is strictly inferior than $\pi$. To this purpose, we recall the technique that the authors have used in \cite{SV21A}.

With no loss of generality, up to a suitable rotation, we can suppose that the spherical triangle ${{\cal{T}}} = ABC$ has the centroid $C_{{\cal{T}}}=(A+B+C)/\|A+B+C\|_2$ cohinciding with the North Pole $(0,0,1)$ and furthermore is strictly contained in the northern-hemisphere, i.e. it does not contain any point of the equator.

Then, if $f \in C({{\cal{T}}})$, $g(x,y)=\sqrt{1-x^2-y^2}$,
$$I_{{{\cal{T}}}}:=\int_{{\cal{T}}} f(x,y,z)d\sigma =  \int_{{{\cal{T}}}^{\perp}} f(x,y,g(x,y)) \frac{1}{g(x,y)}dxdy, $$
where ${{\cal{T}}}^{\perp}$ is the projection of ${{\cal{T}}}$ onto the xy-plane, that is the curvilinear triangle whose vertices, say $\hat{A}$, $\hat{B}$, $\hat{C}$, are the xy-coordinates of ${A}$, ${B}$, ${C}$, respectively (see Figure \ref{fig1}).
%

{ {
The sides of ${{\cal{T}}}^{\perp}$ are arcs of ellipses centered at the origin, being the projections of great circle arcs. Then we can split the planar integral into the sum of the integrals on three elliptical sectors ${\cal{S}}_i$ with $i=1,2,3$, obtained by joining the origin with the vertices $\hat{A}$, $\hat{B}$, $\hat{C}$, namely
$$I_{{{\cal{T}}}}=\int_{{{\cal{T}}}^{\perp}} f(x,y,g(x,y)) \frac{1}{g(x,y)}dxdy=\sum_{i=1}^3 \int_{S_i} f(x,y,g(x,y)) \frac{1}{g(x,y)}dxdy
$$

If the purpose is to compute an algebraic rule over ${\cal{T}}$, with algebraic degree of exactness $ADE=n$, positive weights and internal nodes (i.e. rules of PI-type) we face the problem that while $f$ is a polynomial of total degree $n$ then being $g(x,y)=\sqrt{1-x^2-y^2}$, we have that in general $f(x,y,g(x,y)) \frac{1}{g(x,y)}$ may not be a polynomial.
}}
%

To see this properly, let 
\begin{itemize}
\item $
f(x,y,z)=x^{\alpha} y^{\beta} z^{\gamma}, \quad 0 \leq {\alpha} + {\beta} + {\gamma} \leq n, \,\, {\alpha}, {\beta}, {\gamma} \in \mathbb{N}
$;
\item $g(x,y)=\sqrt{1-x^2-y^2}$.
\end{itemize}
Thus
$$
f(x,y,g(x,y)) \frac{1}{g(x,y)}=x^{\alpha} y^{\beta} (1-x^2-y^2)^{(\gamma-1)/2}.
$$
In particular, if $\gamma$ is 
\begin{itemize}
\item {\sl{odd}} then $f(x,y,g(x,y)) \frac{1}{g(x,y)}$ is a {\sl{polynomial}} of degree at most $n$,
\item {\sl{even}} then $f(x,y,g(x,y)) \frac{1}{g(x,y)}$ is {\sl{$1/g$ multiplied for a polynomial}} of degree at most $n$.
\end{itemize}

To overcome this problem, having in mind to produce a formula that is numerically exact over bivariate polynomials of a certain total degree, we proceed as follows.

First, we approximate $1/g$ by a polynomial $p_{\epsilon}$ of degree $m=m_{\epsilon}$ such that $|p_{\epsilon}-1/g| \leq \epsilon \cdot 1/|g|$. Next, since $f/g \approx f \cdot p_{\epsilon} \in {\mathbb{P}}_{n+m}$, we integrate $f \cdot p_{\epsilon}$ instead of $f/g$ on the elliptical sectors $S_i$, $i=1,2,3$.

Finally, as observed before,  determining a rule of algebraic degree of exactness $n+m$ over each elliptical sector ${\cal{S}}_i$, $i=1,2,3$, with internal nodes, and positive weights then we have a rule on ${{\cal{T}}}^{\perp}:=\cup_{i=1}^3 S_i$ with nodes $(x_k,y_k)_{k=1,\ldots,M}$, weights $w_{k=1,\ldots,M}$ of PI-type.


Mapping back the nodes on the sphere, we have a rule over the spherical triangle that is near algebraic with ADE $n$ since
{\sl{
\begin{equation}
\int_{\cal T} f(x,y,z) d\sigma \approx \sum_{j=1}^{M} \frac{w_i}{\sqrt{1-x_j^2-y_j^2}} f(x_j,y_j,\sqrt{1-x_j^2-y_j^2}).
\end{equation}
}}

At this point, there are some aspects of this procedure that require some further explanation (for details, see also {\cite{SV21A}}).

About the computation of a cubature formula with prescribed degree of exactness over each elliptical sector ${\cal{S}}_i$, since each ${\cal{S}}_i$ is an affine transformation of a circular sector ${{\cal{S}}}^*_i$ of the unit-disk, it is sufficient to obtain  a formula on ${{\cal{S}}}^*_i$ as described in {\cite{DFM}} and map it to ${{\cal{S}}}_i$ (some care on the weights that must be multiplied by absolute value of the transformation matrix determinant). 

\begin{figure}[!ht]
\centering
{\includegraphics[scale=0.4,clip]{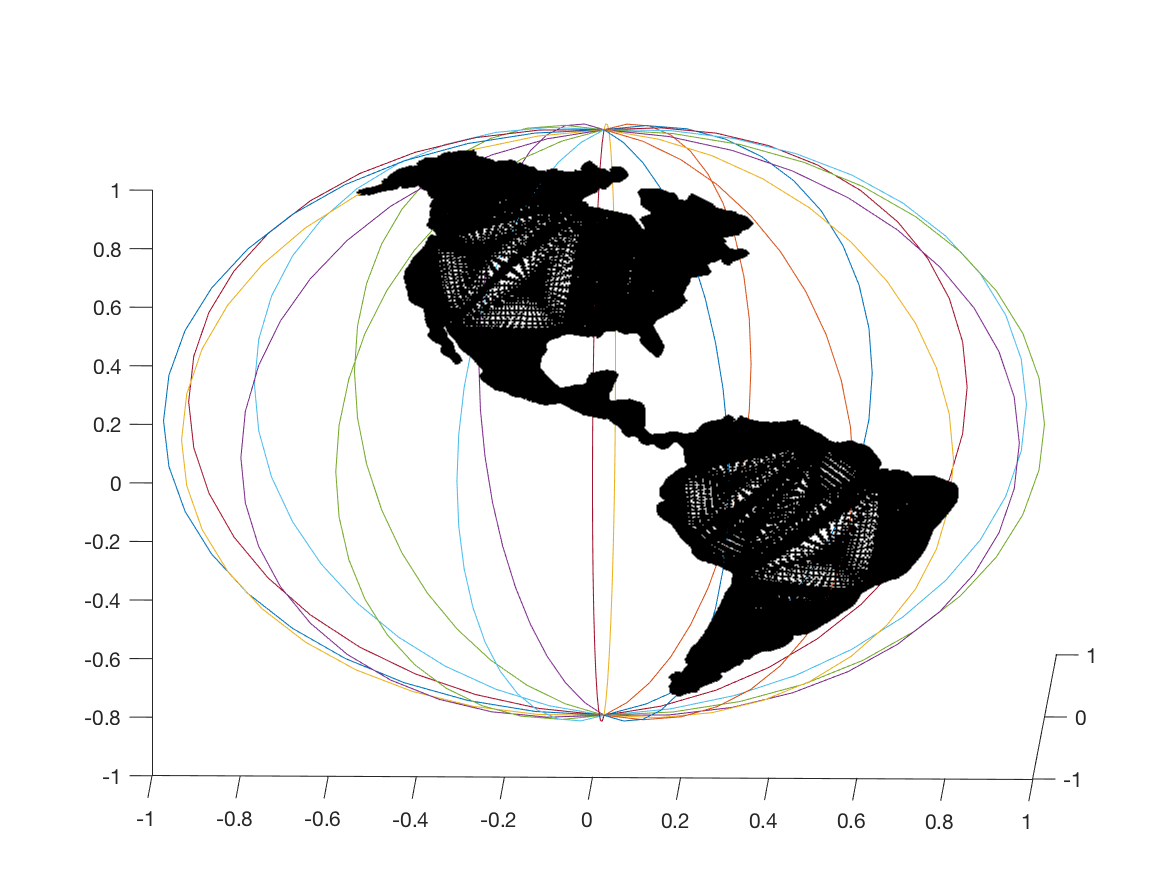}}
 \caption{The 380544 nodes of a rule of PI-type with ADE=8, on a spherical polygon approximating South and North America, before Caratheodory-Tchakaloff compression. The computation of this formula requires about 10 seconds.}
\label{fig2}
\end{figure}

\begin{figure}[!ht]
\centering
{\includegraphics[scale=0.4,clip]{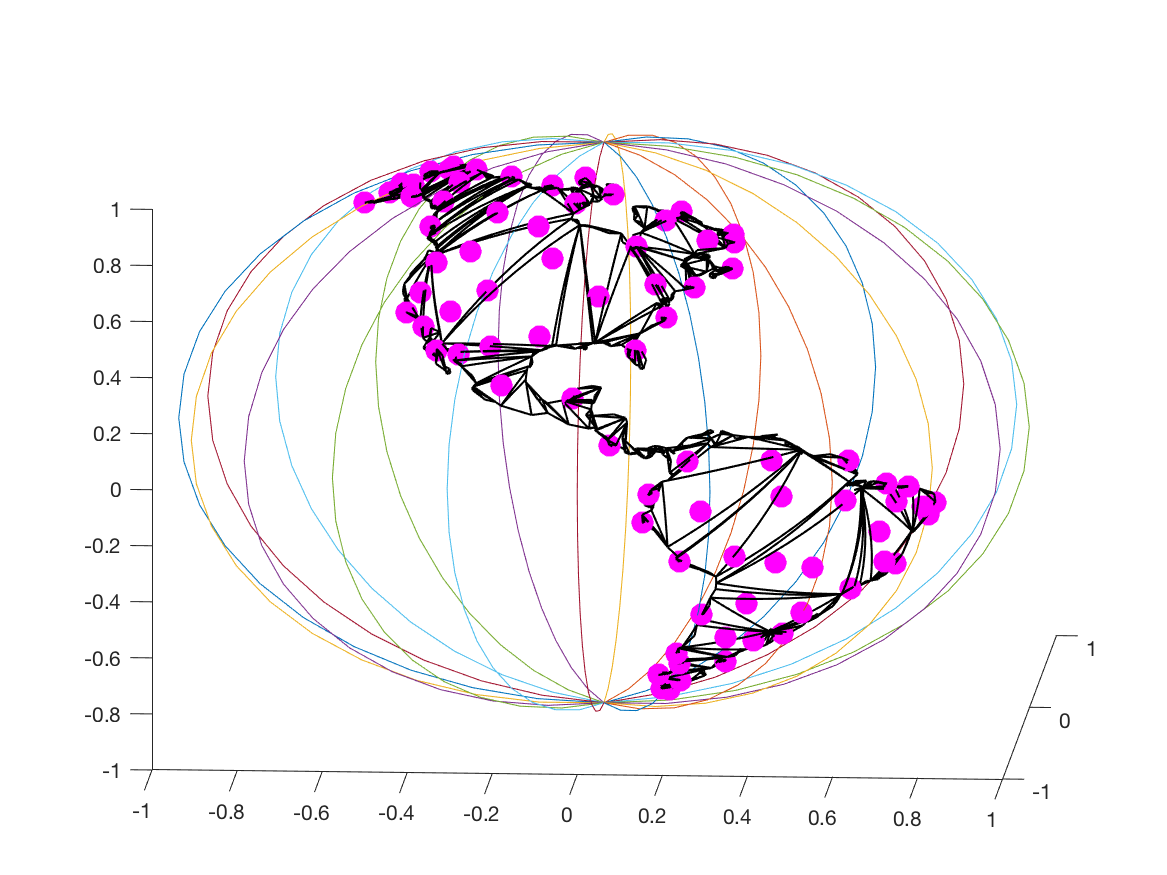}}
 \caption{Triangulation of a spherical polygon (967 spherical triangles) and cubature rule of PI-type with ADE=8, 81 points, after Caratheodory-Tchakaloff compression (magenta). The compression cputime requires approximatively 3.5 seconds.}
\label{fig3}
\end{figure}

Another issue is the computation of $m=m_{\epsilon}$. To this purpose, setting $\epsilon = 10^{-15}$, it turns out that it is sufficient to find what is a (small) degree $m_{\epsilon}$ for which there exist $p_{m_{\epsilon}} \in {\mathbb{P}}_m$ that approximates  
$
{\sl{1/\sqrt{1-t}}}$ where $t \in [0,\max\{\|{\hat{A}}\|_2^2,\|{\hat{B}}\|_2^2,\|{\hat{C}}\|_2^2\}]
$, with a relative error at most $\epsilon$ in $\infty$-norm. Such $m=m_{\epsilon}$ can be quickly estimated by Chebfun {\cite{CHF}} and stored in tables.

We end this section by recalling some basics on the Caratheodory-Tchakaloff rule compression. From the nodes $\{P_k\}_{k=1,\ldots,M}$ and weights $\{w_k\}_{k=1,\ldots,M}$ of the PI-type rule on ${\cal{P}}$, with $M \geq (n+1)^2$, we can extract one with the same features but with cardinality at most $(n+1)^2$, see {\cite{SV15}}. This is guaranteed by Tchakaloff theorem on positive cubature, that in the framework of discrete measures can be demonstrated via the Caratheodory theorem on conical linear combinations of finite-dimensional vectors (applied to the columns of the Vandermonde-like matrix of a moment matching system) and one of these solutions can be obtained numerically by means of Lawson-Hanson algorithm {\cite{LH95}}. We observe that there are several Matlab implementations of the latter procedure, as the Matlab built-in routine {\tt{lsqnonneg}} or that proposed by Slawski in {\cite{SL}}. A new recent fast variant has been introduced in {\cite{DMV}}, and view of its performance, we will use it in our cubature routines to obtain one of these so called {\sl{Caratheodory-Tchakaloff}} formula.

In order to show our approach to a real world problem,  in Figure {\ref{fig2}}, we consider the spherical polygon ${\cal{P}}$ that is an approximation of South and North America (without considering minor islands, rivers and lakes), we plot the 380544 cubature rule nodes of a formula with a numerical ADE equal to 8,  and then in Figure {\ref{fig3}} we exhibit the spherical triangulation and the extracted Caratheodory-Tchakaloff rule with 81 points (exactly the dimension of the polynomial space of total degree $8$ on ${\mathbb{S}}_2$) and again numerical ADE equal to 8.

\subsection{Numerical experiments}

We have implemented the procedure described above in Matlab, so determining the nodes and the weights of a cubature rule over the spherical polygon ${\cal{P}}$ that approximatively integrates all the polynomials $p \in {\mathbb{P}}^3_n({\mathbb{S}}_2)$, for a fixed $n \in {\mathbb{N}}$. The codes are freely available at \cite{SH}. 

A first problem is that there is no known way to achieve $I_{{\cal{P}}}(f)=\int_{{\cal{P}}}f(x,y,z) d\sigma$, when ${\cal{P}}$ has such a complicate geometry. To overcome this problem we have also implemented a basic adaptive code, that computes $I_{{\cal{P}}}(f)$ with a relative error of $10^{-14}$.

As region, we consider Australia (without Tasmania as well as minor islands), whose boundary is provided by the Matlab Mapping toolbox, via longitude-latitude coordinates. Though our software can manage domains that are not simply connected, we did not consider, for sake of simplicity, lakes and rivers. This complex spherical polygon consists of 169 vertices and is the union of 167 spherical triangles whose interiors do not overlap. 

Though the Earth radius is not unitary, in our battery of tests we suppose that ${{\cal{P}}} \subset {\mathbb{S}}_2$. We observe that it is not restrictive, since an integral on a sphere of radius $r$ can be easily reformulated as one on ${\mathbb{S}}_2$.

\begin{figure}[!ht]
\centering
{\includegraphics[scale=0.3,clip]{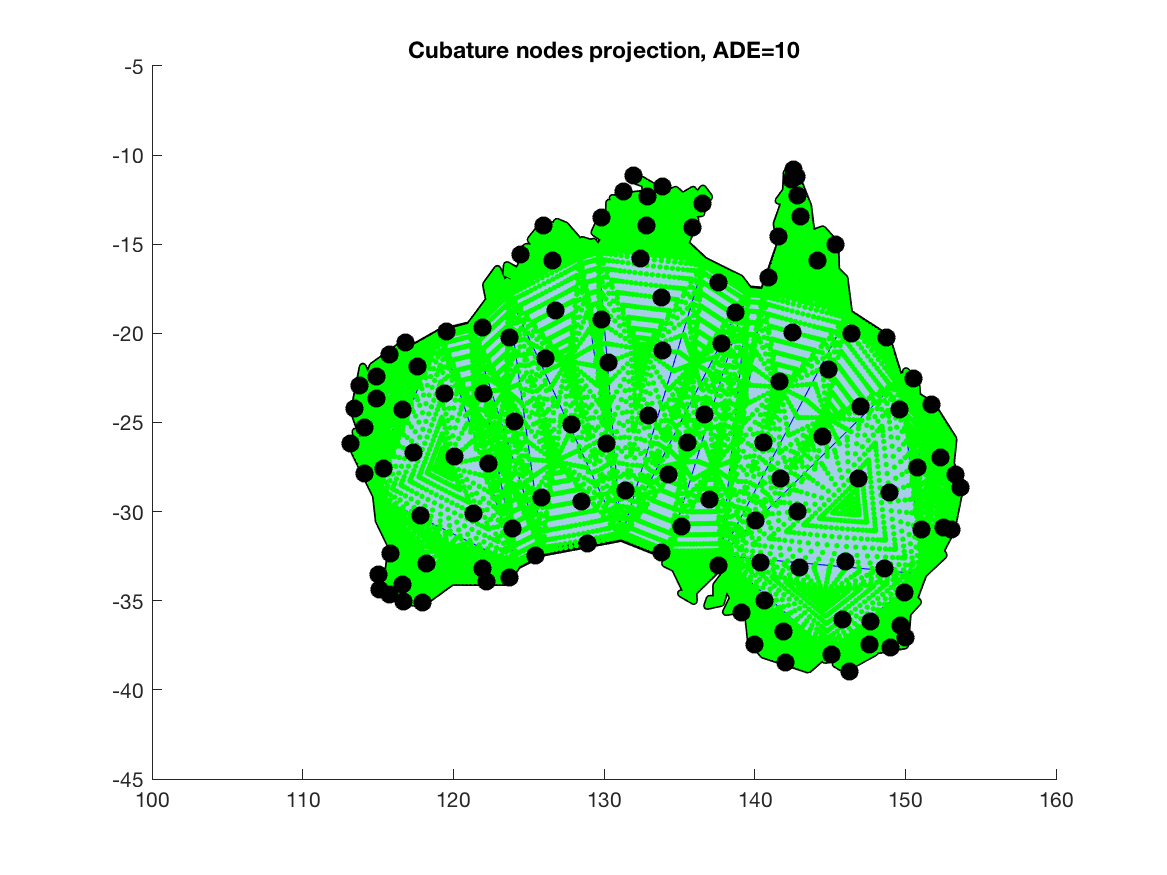}}
 \caption{cubature nodes of a formula of PI-type having ADE equal to 10, on a coarse approximation of Australia using 169 points of ${\mathbb{S}}_2$. To this purpose, we determine a triangulation of such a spherical polygon in 167 spherical triangles, obtained a first rule of PI-type with 82413 nodes, and then a compressed one with 121, with moments error of $\approx 5\cdot 10^{-15}$. The computation of the Caratheodory-Tchakaloff formula took takes 6 seconds.}
 \label{fig5}
 \end{figure}

As for the integrands, setting $(x_0,y_0,z_0)\approx (-0.6325,0.6668,-0.3908) \in {\mathbb{S}}_2$ as an approximation of australian centroid, and $$h(x,y,z) = (x-x_0)^2+(y-y_0)^2+(z-z_0)^2$$ we consider the test functions

\begin{enumerate}
\item $f_1(x,y,z)=1+x+y^2+x^2 \, y+x^4+y^5+x^2 \, y^2 \, z^2$, i.e. a polynomial, of total degree 6;
\item $f_2(x,y,z)=\cos(10\cdot (x+y+z))$, a function with some oscillations;
\item $f_3(x,y,z)=\sin(-h(x,y,z))$, i.e. an analytic function;
\item $f_4(x,y,z)=\exp(-h(x,y,z))$, i.e. an analytic function;
\item $f_5(x,y,z)=((x-x_0)^2 + (y-y_0)^2 + (z-z_0)^2)^{3/2}$, i.e. a $C^{1}({{\cal{P}}})$ function;
\item $f_6(x,y,z)=((x-x_0)^2 + (y-y_0)^2 + (z-z_0)^2)^{5/2}$, i.e. a $C^{2}({{\cal{P}}})$ function.
\end{enumerate}

\begin{figure}
\centering
{\includegraphics[scale=0.4,clip]{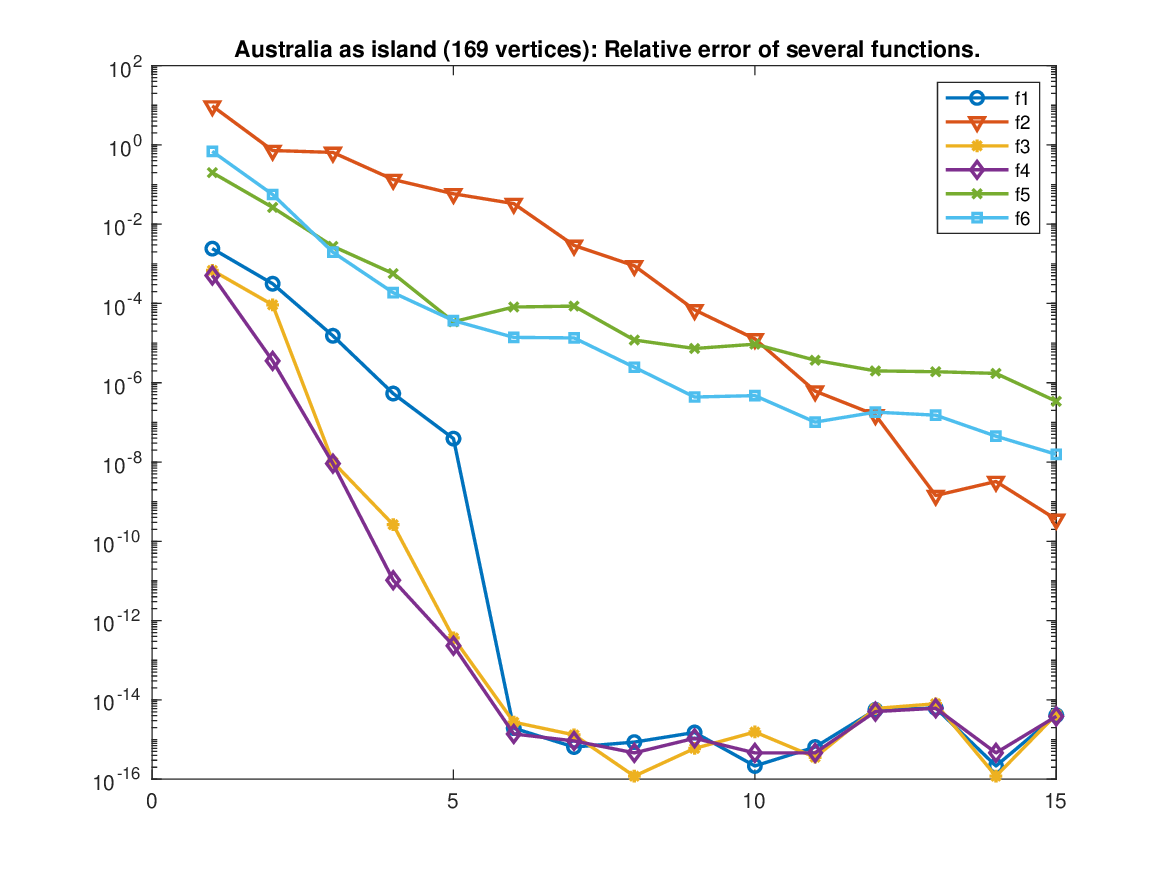}}
 \caption{cubature relative errors, in semilogarithmic scale, obtained by approximating $I_{{\cal{P}}}(f_k)$, for $k=1,2,\ldots,6$, by rules with ADE equal to  $1,2,\ldots,16$.}
 \label{fig6}
 \end{figure}

Due to their properties, we expect that the integration of $f_2$, $f_5$ and $f_6$ by (numerically) algebraic rules will provide {\sl{inferior}}  results w.r.t. $f_1$, $f_3$ and $f_4$. Furthermore, $I_{{\cal{P}}}(f_1)$ will be approximated with a relative error for formulas with degree of exactness at least $6$ of about $10^{-14}$ or less.
The numerical results of this battery of tests are described in Figure \ref{fig6}, and as expected the better results are obtained by $f_1$, $f_3$ and $f_4$ while the convergence is slower for $f_2$, $f_5$ and $f_6$.

Finally, in Table {\ref{T1}}, we have listed the cardinalities of the rules before and after the Caratheodory-Tchakaloff compression, their ratios, and the cputime required by the numerical procedure to determine the Caratheodory-Tchakaloff rule.

These tests show that the rule compression is remarkable and that for mild ADE the cputime necessary to produce a formula of PI type with low cardinality is acceptable, especially in view of the complicated geometry of the integration domain.

All our tests have been performed on a 2.7 GHz Intel Core i5 with 16 GB of RAM, using Matlab 2022a, in which is available the {\tt{polyshape}} environment.
 
 \begin{table}
\begin{center}
\footnotesize
\begin{tabular}{| c | c | c | c | c || c | c | c | c | c |}
\hline
$n$ & \# basic & \# comp & Cratio & CPU & $n$ & \# basic & \# comp & Cratio & CPU  \\
\hline
$1$ & 20613 & 4 & 5153:1 & 0.6s & $9$ & 72441 & 100 & 724:1 & 3.9s \\
$2$ & 25965 & 9 & 2885:1 & 0.6s & $10$ & 82413 & 121 & 681:1 & 6.1s \\
$3$ & 30564 & 16 & 1910:1 & 0.7s & $11$ & 90408 & 144 & 627:1 & 14s \\
$4$ & 37071 & 25 & 1482:1 & 0.8s & $12$ & 101535 & 169 & 600:1 & 17s \\
$5$ & 42519 & 36 & 1181:1 & 1.1s & $13$ & 110379 & 196 & 563:1 & 22s \\
$6$ & 50181 & 49 & 1024:1 & 1.3s & $14$ & 122661 & 225 & 545:1 & 33s \\
$7$ & 56478 & 64 & 882:1 & 1.8s & $15$ & 132354 & 256 & 517:1 & 48s \\
$8$ & 65295 & 81 & 806:1 & 2.6s & $16$ & 145791 & 289 & 504:1 & 62s \\
\hline
\end{tabular}
\caption{\small{Cardinalities of the rule before and after the compression, respectively $\#${\tt{basic}} and $\#${\tt{comp}}, compression ratio {\tt{Cratio}} and CPU time in seconds.}}
\label{T1}
\end{center}
\end{table}

\section {Hyperinterpolation on spherical polygons}
Hyperinterpolation has been introduced in the seminal paper by I.H. Sloan in 1995 (see \cite{SL95}), in the framework of approximation of multivariate functions, and is essentially equivalent to a truncated Fourier expansion in a series of orthogonal polynomials for some discrete or continuous measure.

From then, many aspects has been deepened, either theoretically either from the implementative point of view, proposing hyperinterpolation as a valuable alternative to interpolation, without the issue of determining for the latter a good set of points, e.g. with low Lebesgue constant. In particular, the application to many 2D and 3D domains has been explored, such as balls, cubes, unit sphere, but also less standard ones as polygons, circular sections {\cite{SV17}}, spherical rectangles {\cite{SV18}} and spherical triangles {\cite{SV21B}}.

Denoting by ${\mathbb{P}}^d_n(\Omega)$ the subspace of $d$-variate polynomials of total-degree not exceeding $n$, restricted to a compact set or manifold $\Omega \subset R^d$ {\sl{w.r.t. a given measure}} $d \mu$ on $\Omega$, given 
\begin{itemize}
\item an {\sl{orthonormal basis of ${\mathbb{P}}^d_n(\Omega)$}}, say $\{p_j\}$, $1 \leq j \leq N_n = {\mbox{ dim}} ({\mathbb{P}}^d_n (\Omega))$,
\item a {\sl{cubature formula}} exact for ${\mathbb{P}}_{{2n}}^d (\Omega)$ with nodes $X = \{x_i\} \subset \Omega$ and positive weights ${\bf{w}} = \{w_i\}$, $1 \leq i \leq M$ with $M \geq  N_n$, 
\item the scalar product $ \langle f,g\rangle_M = \sum_{i=1,\ldots,M} w_i f(x_i) g(x_i)$,
\end{itemize}
the {\it{hyperinterpolation}} of $f \in C(\Omega)$ is
\begin{eqnarray}\label{H1}
({\cal{L}}_n f)(x) &=& \sum_{j=1}^{N_n} \langle f,p_j \rangle_M p_j(x) =  {\sl{\sum_{i=1}^{N_n}  w_i f(x_i) \sum_{j=1}^{N_n} p_j(x_i) p_j(x)}}. \nonumber
\end{eqnarray}

Applying this notion to our spherical setting,  if $\Omega \equiv {\cal{P}}$ is spherical polygon strictly contained in a hemisphere, $d \mu=d \sigma$ is the surface measure on the sphere, then
\begin{enumerate}
\item by means of the routines in our Matlab package {\tt{dCATCH}}, freely available at \cite{VdC}, we determine the required {\sl{orthonormal basis}} $\{p_j\}$,
\item next we compute a cubature rule of degree {\sl{$2n$}} on $\Omega$ as previously shown,
\item finally we get the hyperinterpolant of degree {\sl{$n$}}.
\end{enumerate}

Items one and three require some explanation.  Since the cubature formula with nodes $\{x_i\}_{i=1,\ldots,M}$ and weights $\{w_i\}_{i=1,\ldots,M}$ has ADE equal to $2n$, if $\{p_j\}_{j=1,\ldots,N_n}$ is orthonormal w.r.t. this discrete measure, being exact for all the polynomials of total degree at most $2n$, it turns out that it is orthonormal also w.r.t. $d\sigma$. Next, denoting by $V_n \in {\mathbb{R}}^{M \times N_n}$ the Vandermonde matrix of the spherical harmonics basis $\{\phi_i\}_{i=1,\ldots,N_n}$ evaluated at $\{x_i\}_{i=1,\ldots,M}$, we first determine the QR factorization of ${\sqrt{W}} V_n$, i.e. a unitary matrix $Q \in  {\mathbb{R}}^{M \times N_n}$ and a upper triangular one $R \in {\mathbb{R}}^{N_n \times N_n}$ such that ${\sqrt{W}} V_n=QR$ and then set 
$$
(p_1,\ldots,p_{N_n})=(\phi_1,\ldots,\phi_{N_n})R^{-1}.
$$
Since $U_n=V_n R^{-1}=(p_j(x_i))$ is the Vandermonde matrix w.r.t. the basis $\{p_k\}_{k=1,\ldots,N_n}$, it can be easily checked that  the latter is actually orthonormal w.r.t. the scalar product defined by the discrete measure generated by the cubature formula.
Observe that $R$ is non singular since $V_n$ is full rank in view of the fact that the set $X = \{x_i\} \subset \Omega$ is ${\mathbb{P}}_n(\Omega)$-determining.

Furthermore, the vector of hyperinterpolation coefficients ${\boldsymbol{c}}^{\ast}=(c^{\ast}_j)_{j=1,\ldots,N_n}$, such that $({\cal{L}}_n f)(x)=\sum_{j=1}^{N_n} c^{\ast}_j p_j(x)$, that can be regarded as Fourier coefficients, can be conveniently computed in matrix form, being
$$
{\boldsymbol{c}}^*=U^T_nWf=({\sqrt{W}} U_n)^T {\sqrt{W}}{\bf{f}}=Q^T {\sqrt{W}} {\boldsymbol{f}}, \quad {\boldsymbol{f}}=(f(x_i))_{i=1,\ldots,M}.
$$

Concerning the quality of the approximation, there are many theoretical aspects that deserve to be cited, and in particular that if $f$ is a continuous function in $\Omega$ then
\begin{equation}\label{hi1}
\|f-{\cal{L}}_nf\|_{L^2(\Omega)} \leq 2 \sqrt{\sigma(\Omega)} E_n(f,\Omega), \quad E_n(f,\Omega)=\min_{p \in \mathbb{P}_n(\Omega)} \|f-p\|_{\infty},
\end{equation}
so implying that if $p \in {\mathbb{P}}^d_n(\Omega)$ then ${\cal{L}}_n p=p$, as well as
\begin{equation}\label{hi2}
\|f-{\cal{L}}_n\|_{L^{\infty}(\Omega)} 	\leq (1+\|{\cal{L}}_n\|)E_n(f,\Omega)
\end{equation}
where 
$$
\|{\cal{L}}_n\|=\sup_{f \neq 0} \frac { \|{\cal{L}}_n f\|_{L^{\infty}(\Omega)}}  { \|f\|_{L^{\infty}(\Omega)}} .
$$ 
For additional more technical details, see {\cite{SV17}}, {\cite{SV21B}}.

{\vspace{0.2cm}}

Recently some hyperinterpolation variants have been studied, having in mind between various scopes also to improve the approximation quality in case of noisy data. We make a brief glance on some of them, as the filtered, Lasso and hybrid variants, useful later to comprehend the comparisons made in the numerical section. 

In {\it{filtered hyperinterpolation}} {\cite{SLW12}}, one introduces a {\it{filter function}} $h \in \mathcal{C}([0,+\infty))$ that satisfies
\begin{equation*}
h(x)=\begin{cases}
1, &\text{ for }x\in[0,1/2],\\
0, &\text{ for }x\in[1,\infty).
\end{cases}
\end{equation*}
It is straightforward to notice that depending on the behaviour in $[1/2,1]$, one can define many filters, e.g. as used in {\cite{ARS}}
\begin{equation}\label{equ:filtered}
 h(x) = \left\{\begin{array}{cc}
1, & x \in [0,\frac{1}{2}] ,\\
\sin^2(\pi x), & x \in [\frac{1}{2},1]\\
0, &\text{ for }x\in[1,\infty).
\end{array}
\right.
\end{equation}

Denoting by $\lfloor \cdot \rfloor$ the floor function, by $h$ the choosen {\it{filter}} and by $\langle f,g \rangle_M$ a discrete scalar product defined by an $M$-point cubature rule of PI-type in $\Omega$ with algebraic degree of exactness  $L-1+\lfloor L/2 \rfloor$, the {\it{filtered hyperinterpolant}} ${\cal{F}}_{n} f \in \mathbb{P}_{n-1}(\Omega)$ of $f$ consists in
\begin{equation}\label{filt}
{\cal{F}}_{n} f:=\sum_{j=1}^{N_n} h\left({\frac{\deg p_{j}}{n}}\right) \langle{f,p_{j}}\rangle_M \, \, p_{j}.
\end{equation}
It can be easily seen, in view of the definition of the filter, that if $p \in {\mathbb{P}}_{\lfloor n/2 \rfloor}^d(\Omega)$ then ${\cal{F}}_{n} p=p$. Next, in \cite{L21} the authors show that (distributed) filtered hyperinterpolation can reduce weak noise.

One of the potential difficulties in (\ref{filt}) is the evaluation of $\deg p_{j}$, $j=1,\dots,N_n$.
This is straightforward for the classical spherical harmonics basis on ${\mathbb{S}}^2$, since it is a {\it{triangular basis}}, i.e. for $k=0,\ldots,d$, $j=k^2+1,\dots,(k+1)^2$, we have $\deg p_{j}=k$. 

In the case of the orthonormal basis $\{p_j\}_{j=1,\ldots,N_n}$ on a spherical polygon ${\cal{P}}$, one has to take into account that it is computed numerically as specified above, by means of QR factorisation, starting from the Vandermonde matrix defined by spherical harmonics. Anyway, since ${R}^{-1}$ is an upper triangular non-singular matrix, then the so determined orthonormal basis $(p_1,\ldots,p_{N_n})$ of ${\mathbb{P}}_n^3({\cal{P}})$ is also triangular and again for $k=0,\ldots,d$, $j=k^2+1,\dots,(k+1)^2$, we have that $\deg p_{j}=k$.

An alternative approach consists in {\it{Lasso hyperinterpolation}} \cite{an2021Lasso} that denoises perturbed data and attempts to dismiss less relevant discrete Fourier coefficients.

To this purpose, let us define the {\it{soft thresholding operator}} as
\[\softoperatordef{a}:=\max(0,a-k)+\min(0,a+k),\]
where $k \geq 0$, 
that is
\[
\softoperatordef{a}=\left\{
\begin{array}{cl}
a+k, & {\mbox{ if }} a < -k ,\\
0,  & {\mbox{ if }} -k \leq a \leq k, \\
a-k, & {\mbox{ if }}  a > k. \\
\end{array}
\right.
\]
and suppose that a discrete scalar product $\langle f,g \rangle_M$ is given by an $M$-point cubature rule of PI-type in $\Omega$ with algebraic degree of exactness $2n$.

Then the \emph{Lasso hyperinterpolation} of $f$ onto $\mathbb{P}_{L}(\Omega)$ is defined as
\begin{equation}\label{Lasso}
{\cal{L}}_L^{\lambda}{f}:=\sum_{j=1}^{N_n} {\cal{S}}_{\lambda \mu_{j}} (\langle{f,p_{j}}\rangle_M) \, \, p_{j}.
\end{equation}
where $\lambda>0$ is the regularization parameter and $\{\mu_{k}\}_{k=1}^{N_n}$ is a set of positive penalty parameters.

By the definition of $\softoperatordef{a}$, if $|\langle f,p_k \rangle_M| \leq \lambda \mu_{k}$ then ${\cal{S}}_{\lambda \mu_{k}} \discreteinner{f,p_k}=0$, so explaining how it tries to dismiss less relevant discrete Fourier coefficients, slightly modifying the remaining ones.
Differently from the classical hyperinterpolation and filtered hyperinterinterpolation, the Lasso operator is not in general a projection to a certain polynomial space and is not invariant under a change of basis.

In {\cite{an2021Lasso}} the authors have shown theoretically and numerically the effectiveness  ${\cal{L}}^{\lambda}_n$ in the case of noisy data, taking advantage of the connection between ${\cal{L}}_n^{\lambda}{f}$ and a certain $\ell_1-$regularized least squares problem.

In order to combine the features of filtered and Lasso hyperinterpolation, in {\cite{ARS}} the authors introduced the so called {\it{hybrid hyperinterpolation}}.

Suppose that the discrete scalar product $\langle f,g \rangle_M$ is determined as above by an $M$-points quadrature rule of PI-type in $\Omega$ with algebraic degree of exactness $2n$. The \emph{hybrid hyperinterpolation} of $f$ onto $\mathbb{P}_{n}^d (\Omega)$ is defined as
\begin{equation}\label{hyb}
\hybridhyper{f}:=\sum_{j=1}^{N_n} h\left(\frac{\deg p_{j}}{n} \right)\mathcal{S}_{\lambda \mu_{j}} (\discreteinner{f, p_{j}})p_{j}, \quad \lambda>0,
\end{equation}
where $h(\cdot)$ is a filter function, $\{\mathcal{S}_{\lambda \mu_{j}}(\cdot)\}_{j=1}^{N_n}$ are soft thresholding operators in which $\mu_k >0$, $k=1,\ldots,N_n$.

Similarly to the case of the Lasso hyperinterpolation, by means of a connection between $\hybridhyper{f}$ and a certain ${\ell^2_2}-\ell_1$regularized least squares problem, in {\cite[Thm. 4.2]{ARS}} it is shown how hybrid hyperinterpolation acts on noisy data, with advantages in term of sparsity of its polynomial coefficients.  

In the numerical section we report some examples, in which we show the performance of the classical hyperinterpolation as well as a comparison in problems dealing with noisy data.

\begin{Remark}
In this section we have not mentioned for sake of brevity some relevant hyperinterpolation variants as {\it{Hard Thresholding}} {\cite{AHT}} or Tikhonov regularized least squares approximation that continuously shrinks all hyperinterpolation coefficients without dismissing less relevant ones. 
\end{Remark}

\begin{Remark}
We observe that, though ${\cal{L}}_n$ requires the availability of on algebraic cubature rule with degree of precision $ADE=2n$, some recent papers are considering the case when instead we have that $ADE< 2n$ and how this choice affects the approximation quality (see, e.g., \cite{AWBIT},\cite{AW}).
\end{Remark}

\subsection{Numerical experiments}

In our first numerical examples, we consider as spherical polygon $\cal P$ that obtained as course map of Australia (without taking into account its smaller islands, lakes and rivers), scaled in $S_2$ and as functions those $f_k$, $k=1,\ldots,6$ defined above.

\begin{figure}[ht!]
\centering
{\includegraphics[scale=0.4,clip]{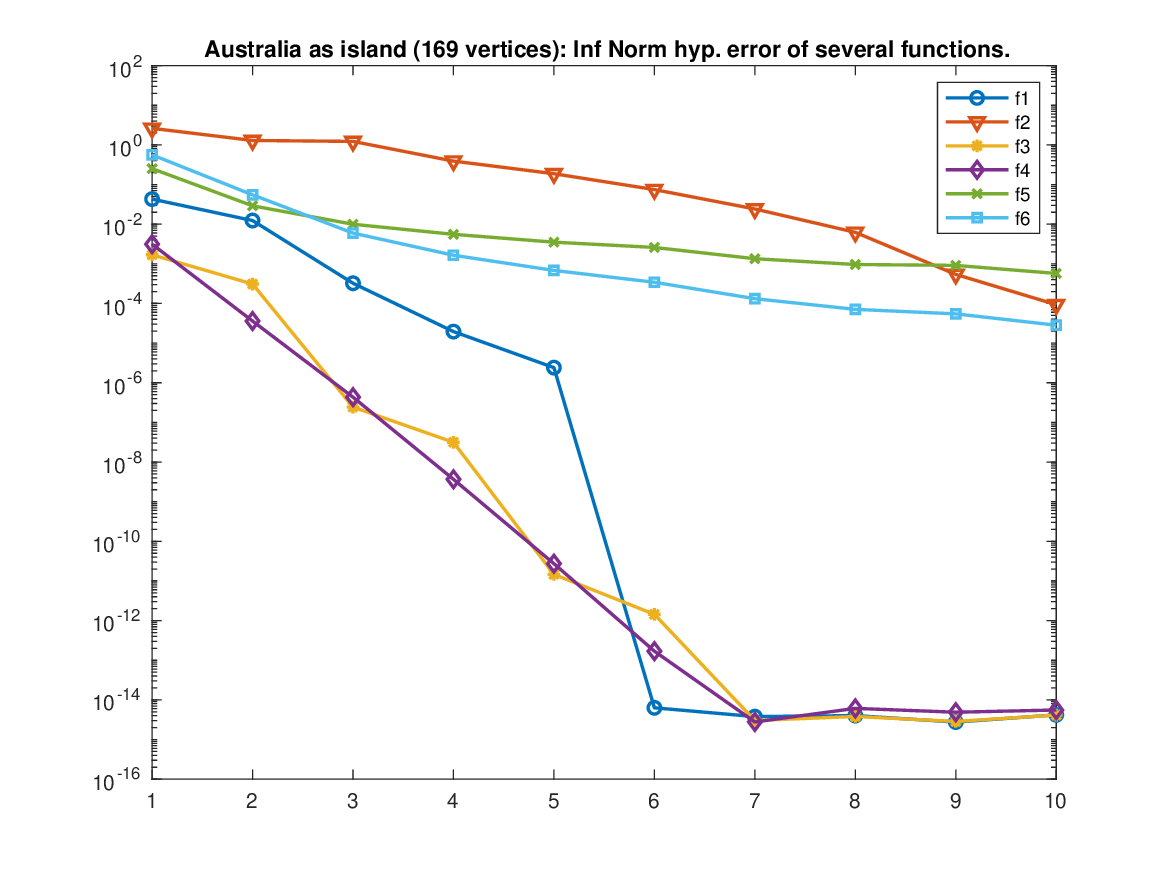}}
 \caption{Inf-Norm hyperinterpolation relative error of several test functions over a coarse approximation of Australia.}
 \label{hL2a}
 \end{figure}

In  Figure {\ref{hL2a}} we report an approximation of the relative error
$$
\frac{\|f_k(X)-({\cal{L}}_n f_k)(X)\|_{\infty}}{\|f_k(X)\|_{\infty}} \approx \frac{\|f_k-{\cal{L}}_n f_k\|_{\infty}}{\|f_k\|_{\infty}}
$$
where $n=1,\ldots,10$ and $X$ is the set of 82413 nodes of the rule on ${{\cal{P}}}$, with ADE equal to $10$.
 
As somehow expected, the quality of the approximation of $f_2$, $f_5$ and $f_6$ via hyperinterpolation is {\sl{inferior}} to that of $f_1$, $f_3$ and $f_4$. Furthermore, since $f_1 \in {\mathbb{P}}_6$, in exact arithmetic we have ${\cal{L}}_n f_1 \equiv f_1$ for $n \geq 6$, so explaining why we get in this case relative errors close to $10^{-15}$.  
 
\begin{figure}[ht!]
\centering
{\includegraphics[scale=0.4,clip]{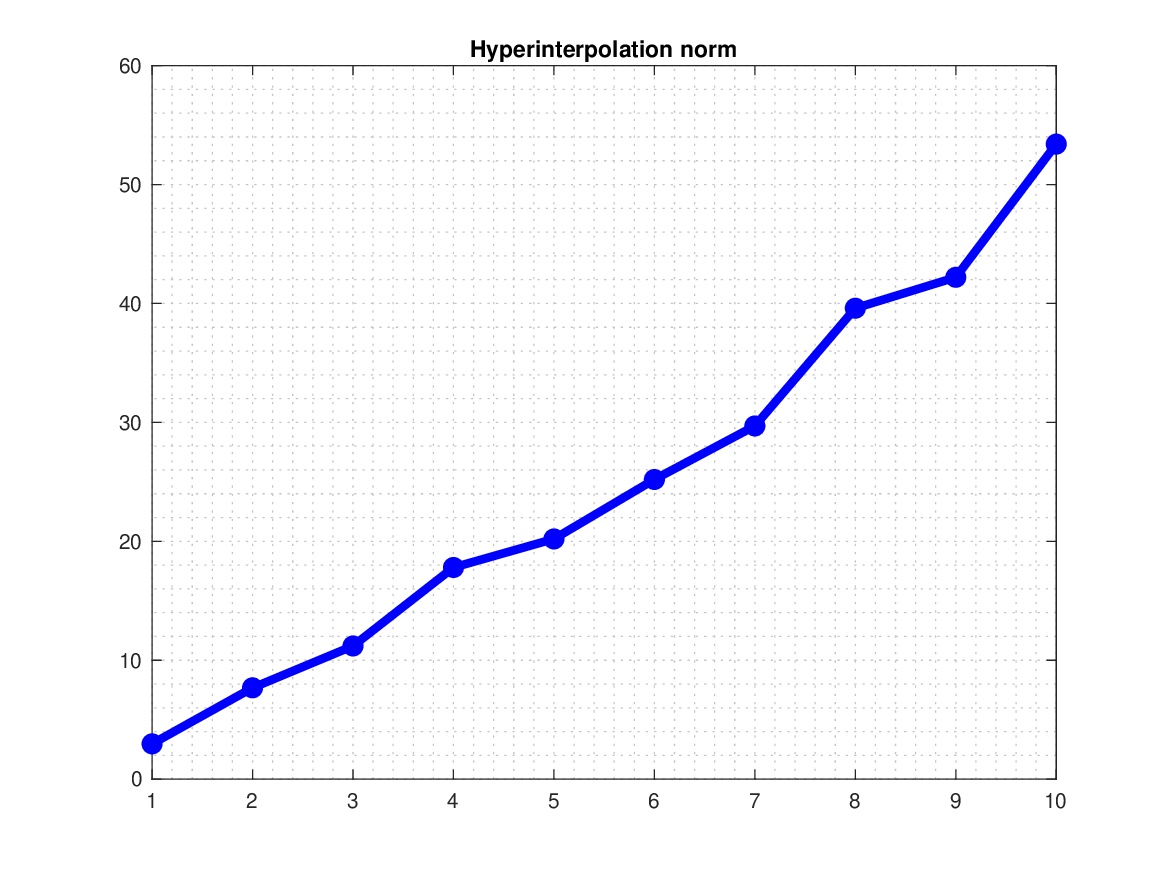}}
 \caption{Approximation of the hyperinterpolation operator norm $\|{\cal{L}}_n\|$ when ${\cal{P}}$ is a coarse approximation of Australia.}
 \label{hL2b}
 \end{figure}

Next we plot in Figure {\ref{hL2b}} an approximation in semilogarithmic scale of $\|{\cal{L}}_n\|_{\infty}$ for $n=1,\ldots,10$, a quantity that as seen above is fundamental to determine, in view of (\ref{hi2}), an upper bound of $\|f-{\cal{L}}_n\|_{\infty} $ w.r.t. the best approximation error $E_n(f,{{\cal{P}}})$, where $f \in C({{\cal{P}}})$.

To this purpose, defining the reproducing kernel $K_n(x,y)=\sum_{i=1}^{N_n} p_i(x) p_i(y)$ and $g_i(x)=w_i K_n(x_i,x)$, it can be proven that 
$$
\|{\cal{L}}_n\|_{\infty}=\max_{x \in {{\cal{P}}}} \sum_{i=1}^M |g_i(x)| \approx \max_{x \in X} \sum_{i=1}^M |g_i(x)|
$$
where $X$ is a sufficiently dense set of points. In our tests, we have choosen again as $X$ the set of 82413 nodes of the rule on ${{\cal{P}}}$, with ADE equal to $10$. The results show that for $n \leq 10$, the operator norm $\|{\cal{L}}_n\|$ is relatively small and thus, by (\ref{hi2}), the approximation error is not too large w.r.t. that of the best approximant $E_n(f,{{\cal{P}}})$.

We report that, when the cubature formula is available, the computation of the hyperinterpolation coefficients ${\bf{c}}^*=(c_k)_{k=1,\ldots,N_n}$, for $n \leq 10$, requires from $5\cdot 10^{-4}$ to $10^{-2}$ seconds.

{\vspace{0.2cm}}

Next we intend to compare the proposed forms of hyperinterpolation, i.e. the classical one as well the filtered, the Lasso and the hybrid one on the case of noisy data in the spherical polygon ${\cal{P}}$ adopted above.

As for filtered hyperinterpolation, following {\cite{ARS}}, we choose $h: [0,1] \rightarrow [0,1]$ defined as

$$
h(x)=\left\{
\begin{array}{cl}
1, &x \in [0, \frac{1}{2}], \\
\sin^2 (\pi \, x), &x \in [\frac{1}{2},1].
\end{array}
\right.
$$

In the experiments reported in Table \ref{noisy1}, we consider noisy evaluations of $$f(x,y,z) = \exp(x^6 \cdot \cos(y+2z))$$ on the nodes $\{ P_j\}_{j=1,\ldots,441}$ of a compressed cubature rule over ${\cal{P}}$ with ADE equal to $20$.

The perturbation is obtained by means of the sum of {\em{Gaussian noise}} ${\cal{N}}(0,\sigma^2)$ from a normal distribution with mean 0
 and standard deviation $\sigma$, and of {\em{impulse noise}} ${\cal{I}}(a)$ that takes uniformly distributed random values in $[-a,a]$ with probability $1/2$.

\noindent
In these tests, we have set $a=0.25$ and $\sigma=0.25$.

Next, being $p_{10}$ the chosen hyperinterpolant of degree $10$ on the perturbed data $\{(P_j,f(P_j)+\epsilon_j)\}_{j=1,\ldots,441}$, we estimate the errors $\|f-p_{10}\|_2$ via a discretization of the $L_2$ scalar product by means of a cubature rule with degree of precision $30$. In Table \ref{noisy1} we report as result the average error over 10 trials.

{\vspace{0.1cm}}

\begin{table}
\begin{center}
  \begin{tabular}{|c| c c c c |}
   \hline
   $$ &  $\lambda \approx 0.011$ &  $\lambda \approx 0.0087$  &  $\lambda \approx 0.0073$  &  $\lambda \approx 0.0061$ \\
    \hline
 	 {\em{filtered}} & $0.070603$  & $    0.067005$  & $  0.069905 $  & $   0.068639$  \\

 	 {\em{Lasso}} &$0.034094 $  & $   0.035766 $  & $   0.037304   $  & $  0.03865$  \\

 	 {\em{hybrid}} & $0.033333  $  & $  0.033969  $  & $  0.034454  $  & $   0.03513$  \\

 	 {\em{hyperint.}} & $0.079629  $  & $   0.076223  $  & $   0.079625  $  & $   0.078205$ \\

 	\hline

 	 {\em{hybrid sparsity}} &$19$ &$29$ &$39$ &$49$ \\
    \hline
  \end{tabular}
  \caption{Average approximation errors in the 2-norm and sparsity of hybrid hyperinterpolation coefficients of noisy versions of $f(x) = \exp((x^6) \cos(y+2z))$, over a spherical polygon {\cal{P}} that approximates the Australia, scaled on the unit-sphere ${\mathbb{S}}_2$, via variants of hyperinterpolation of degree $10$. The perturbation is obtained by the sum of impulse noise with $a=0.25$ and of gaussian noise having standard deviation $\sigma=0.25$.}
  \label{noisy1}
  \end{center}
\end{table}

As for the Lasso and hybrid parameters $\lambda$, we use some of the hyperinterpolant coefficients that are $20$-th, $30$-th, $40$-th and $50$-th in magnitude, while all $\mu_i$ are equal to $1$.

In these tests the performance of Lasso and hybrid hyperinterpolation is superior to that of the classical and filtered hyperinterpolation. Observe that the latter does not depend on the choice of $\lambda$, so explaining why their results remain almost equal varying such parameter. 

The last column shows the sparsity of the hybrid hyperinterpolation coefficients, i.e. the number of those that are non null, varying $\lambda$. As comparison, take into account that filtered and classical hyperinterpolation require $121$ non null coefficients.

As final experiments, we report in Table \ref{noisy2}, the results obtained as in the previous example but fixing $\lambda$ and varying $a$ and $\sigma$, i.e. the amount of noise. In this framework the chosen values of $\lambda$, each time obtained as the $20$-th hyperinterpolation coefficient in magnitude, vary in the interval $[0.0114,0.0150] $.

{\vspace{0.1cm}}

\begin{table}
\begin{center}
  \begin{tabular}{|c| c c c c |}
   \hline
   $$ &  $a=\sigma=0.025$ &  $a=\sigma=0.05$  &  $a=\sigma=0.1$  &   $a=\sigma=0.3$  \\
    \hline
    
     {\em{filtered}} & $        0.00678  $ &$  0.0137  $ &$  0.0276    $ &$ 0.0844$ \\
     {\em{Lasso}} &$  0.0124   $ &$  0.0151   $ &$  0.0200  $ &$   0.0391$ \\
     {\em{hybrid}} & $           0.0124  $ &$   0.0150 $ &$    0.0198  $ &$   0.0376$ \\
    {\em{hyperint.}} & $ 0.0076  $ &$   0.0157  $ &$   0.0313  $ &$   0.0964$ \\
    \hline

 	 {\em{hybrid sparsity}} &$ 19      $ &$        19    $ &$          19        $ &$      19$ \\
    \hline
  \end{tabular}
  \caption{Average approximation errors in the 2-norm and sparsity of hybrid hyperinterpolation coefficients of noisy versions of $f(x) = \exp((x^6) \cos(y+2z))$, over a spherical polygon {\cal{P}} that approximates the Australia, scaled on the unit-sphere ${\mathbb{S}}_2$, via variants of hyperinterpolation of degree $10$. The perturbation is obtained by the sum of impulse noise with several values of $a$ and of gaussian noise having standard deviation $\sigma$, taking as $\lambda$ the $10$-th Fourier coefficient in magnitude of hyperinterpolation.}
  \label{noisy2}
  \end{center}
\end{table}

We observe the advantage of using filtered or classic hyperinterpolation in case of mild noise while it is clear the counterpart when the perturbation is more significative.

We conclude saying that all the Matlab codes for hyperinterpolation on spherical polygons used in this work are available as open-source software at {\cite{SH}}.

%


%

\bibliography{}
\bibliographystyle{plain}

\end{document}